\documentclass[11pt,twoside]{article}

\usepackage{amsbsy,amsfonts,amsmath,amssymb,eucal,mathrsfs}
\usepackage[all]{xy}
\usepackage{pstricks}
\usepackage{pst-plot}
\newtheorem{definition}{Definition}[section]
\newenvironment{defi}{\begin{definition} \rm}{\end{definition}}

\newtheorem{prop}[definition]{Proposition}

\newtheorem{lemm}[definition]{Lemma}
\newtheorem{fact}[definition]{Fact}
\newtheorem{coro}[definition]{Corollary}
\newtheorem{theo}{Theorem}
\newtheorem{notation}[definition]{Notation}

\newtheorem{construction}[definition]{Construction}

\newtheorem{remark}[definition]{Remark}
\newenvironment{rema}{\begin{remark} \rm}{\end{remark}}
\newtheorem{remarks}[definition]{Remarks}

\newtheorem{example}[definition]{Example}

\newtheorem{examples}[definition]{Examples}

\newtheorem{nothing}[definition]{$\!\!$}

\newenvironment{proo}{{\flushleft \it Proof.}}{\hfill $\square$ \vspace{2mm}}

\newtheorem{definition*}{Definition}[section]
\newenvironment{defi*}{\begin{definition*} \rm}{\end{definition*}}
\newtheorem{definitions*}[definition*]{Definitions}
\newenvironment{defis*}{\begin{definitions*} \rm}{\end{definitions*}}
\newtheorem{prop*}[definition*]{Proposition}
\newtheorem{lemm*}[definition*]{Lemma}
\newtheorem{coro*}[definition*]{Corollary}
\newtheorem{theo*}[definition*]{Theorem}
\newtheorem{remark*}[definition*]{Remark}
\newenvironment{rema*}{\begin{remark*} \rm}{\end{remark*}}
\newtheorem{remarks*}[definition*]{Remarks}
\newenvironment{remas*}{\begin{remarks*} \rm}{\end{remarks*}}
\newtheorem{example*}[definition*]{Example}
\newenvironment{exam*}{\begin{example*} \rm}{\end{example*}}
\newtheorem{examples*}[definition*]{Examples}
\newenvironment{exams*}{\begin{examples*} \rm}{\end{examples*}}
\newtheorem{nothing*}[definition*]{$\!\!$}
\newenvironment{noth*}{\begin{nothing*} \rm}{\end{nothing*}}

\begin{document}

\def \sca #1#2{\langle #1,#2 \rangle}
\def\pt{\{{\rm pt}\}}
\def\aut{{\rm Aut}}
\def\ra{\rightarrow}
\def\s{\sigma}\def\OO{\mathbb O}\def\PP{\mathbb P}\def\QQ{\mathbb Q}
 \def\CC{\mathbb C} \def\ZZ{\mathbb Z}\def\JO{{\mathcal J}_3(\OO)}
\newcommand{\G}{\mathbb{G}}
\def\proof{\noindent {\it Proof.}\;}
\def\qed{\hfill $\square$}
\def \uh {{\widehat{u}}}
\def \vh {{\widehat{v}}}
\def \fh {{\widehat{f}}}
\def \wh {{\widehat{w}}}
\def \Wh {{{W_{{\rm aff}}}}}
\def \Wt {{\widetilde{W}_{{\rm aff}}}}
\def \Qt {{\widetilde{Q}}}
\def \Waff {{W_{{\rm aff}}}}
\def \Waffm {{W_{{\rm aff}}^-}}
\def \Wpaff {{{(W^P)}_{{\rm aff}}}}
\def \Wtpaff {{{(\widetilde{W}^P)}_{{\rm aff}}}}
\def \Wtaffm {{\widetilde{W}_{{\rm aff}}^-}}
\def \lh {{\widehat{\lambda}}}
\def \pit {{\widetilde{\pi}}}
\def \lt {{{\lambda}}}
\def \xh {{\widehat{x}}}
\def \yh {{\widehat{y}}}

%%%%%%%%% Commandes pe %%%%%%%%%%%%%%

% alg{\~A}{\AA}{!`}bres norm{\~A}{\copyright}es et anneaux usuels

\newcommand{\N}{\mathbb{N}}
\newcommand{\A}{{\mathbb{A}_{\rm Aff}}}
\newcommand{\Ah}{{\mathbb{A}_{\rm Aff}}}
\newcommand{\At}{{\widetilde{\mathbb{A}}_{\rm Aff}}}
\newcommand{\Ht}{{{H}^T_*(\Omega K^{\ad})}}
\renewcommand{\H}{{{H}^T_*(\Omega K)}}
\newcommand{\Ih}{{I_{\rm Aff}}}
\newcommand{\psit}{{\widetilde{\psi}}}
\newcommand{\xit}{{\widetilde{\xi}}}
\newcommand{\Jt}{{\widetilde{J}}}
\newcommand{\Zt}{{\widetilde{Z}}}
\newcommand{\at}{{\widetilde{A}}}
\newcommand{\Z}{\mathbb Z}
\newcommand{\R}{\mathbb{R}}
\newcommand{\Q}{\mathbb{Q}}
\newcommand{\C}{\mathbb{C}}
\renewcommand{\O}{\mathbb{O}}
\newcommand{\F}{\mathbb{F}}
\newcommand{\p}{\mathbb{P}}
\newcommand{\co}{{\cal O}}

\renewcommand{\a}{{\cal A}}
\newcommand{\az}{\a_\Z}
\newcommand{\ak}{\a_k}

\newcommand{\rc}{\R_\C}
\newcommand{\cc}{\C_\C}
\newcommand{\hc}{\H_\C}
\newcommand{\oc}{\O_\C}

\newcommand{\rk}{\R_k}
\newcommand{\ck}{\C_k}
\newcommand{\hk}{\H_k}
\newcommand{\ok}{\O_k}

\newcommand{\rz}{\R_Z}
\newcommand{\cz}{\C_Z}
\newcommand{\hz}{\H_Z}
\newcommand{\oz}{\O_Z}

\newcommand{\RR}{\R_R}
\newcommand{\CR}{\C_R}
\newcommand{\HR}{\H_R}
\newcommand{\OR}{\O_R}

\newcommand{\re}{\mathtt{Re}}

% alg{\~A}{\AA}{!`}bre

\newcommand{\dual}{{\bf v}}
\newcommand{\com}{\mathtt{Com}}
\newcommand{\rg}{\mathtt{rg}}
\newcommand{\scal}[1]{\langle #1 \rangle}

\newcommand{\g}{\mathfrak g}
\newcommand{\h}{\mathfrak h}
\renewcommand{\u}{\mathfrak u}
\newcommand{\n}{\mathfrak n}
\newcommand{\e}{\mathfrak e}
\newcommand{\plie}{\mathfrak p}
\newcommand{\tlie}{\mathfrak t}
\newcommand{\llie}{\mathfrak l}
\newcommand{\q}{\mathfrak q}
\newcommand{\liesl}{\mathfrak {sl}}
\newcommand{\so}{\mathfrak {so}}

\newcommand{\ad}{{\rm ad}}
\newcommand{\jad}{{j^\ad}}
\newcommand{\id}{{\rm id}}

%%%% Caligraphic letters %%%%%%%%%%%%%%

\def\cA{{\cal A}} \def\cC{{\cal C}} \def\cD{{\cal D}} \def\cE{{\cal E}}
\def\cF{{\cal F}} \def\cG{{\cal G}} \def\cH{{\cal H}} \def\cI{{\cal I}}
\def\cK{{\cal K}} \def\cL{{\cal L}} \def\cM{{\cal M}} \def\cN{{\cal N}}
\def\cO{{\cal O}}
\def\cP{{\cal P}} \def\cQ{{\cal Q}} \def\cT{{\cal T}} \def\cU{{\cal U}}
\def\cV{{\cal V}} \def\cX{{\cal X}} \def\cY{{\cal Y}} \def\cZ{{\cal Z}}

\title{Affine symmetries of the equivariant quantum\\
cohomology ring of rational homogeneous spaces}
 \author{P.E. Chaput, L. Manivel, N. Perrin}

\maketitle

\begin{abstract}
Let $X$ be a rational homogeneous space and let $QH^*(X)_{loc}^\times$ be the
group of invertible elements in the small quantum cohomology ring of
$X$ localised in the quantum parameters. We generalise results of
\cite{cmp} and realise explicitly the map $\pi_1({\rm Aut}(X))\to
QH^*(X)_{loc}^\times$ described in \cite{seidel}. We even prove that
this map is an embedding and realise it in the equivariant
quantum cohomology ring $QH^*_T(X)_{loc}^\times$. We give explicit
formulas for the product by these elements.

The proof relies on a generalisation, to a quotient of the
equivariant homology ring of the affine Grassmannian, of a formula
proved by Peter Magyar \cite{Magyar}. It also uses Peterson's
unpublished result \cite{Peterson} --- recently proved by Lam and
Shimozono in \cite{Lam-Shi} --- on the comparison between the
equivariant homology ring of the affine Grassmannian and the
equivariant quantum cohomology ring.
\end{abstract}

 {\def\thefootnote{\relax}
 \footnote{ \hspace{-6.8mm}
 }
 }

\section{Introduction}

For $G$ a semisimple simply connected algebraic group, the center
$Z$ of $G$ has several interpretations. It may be
canonically identified to the fundamental group $\pi_1(G^{\ad})$ of
the adjoint group $G^{\ad}=G/Z$. Another description is
$$Z\simeq P^\vee/Q^\vee$$
where $P^\vee$ and $Q^\vee$ are the coweight and coroot lattices (see
for example \cite{bou}). Nice representatives of this quotient in
$P^\vee$ are given by the opposites of the minuscule fundamental coweights
$(-\varpi_i^\vee)_{i\in I_m}$ (recall that a dominant coweight
$\lambda$ is minuscule if $\sca{\lambda}{\alpha}=0$ or 1 for any positive
root $\alpha$). Here
$I_m$ denotes the subset of the set $I$ of vertices of the Dynkin diagram
of $G$ parametrising the minuscule coweights.

The group $Z$ can also be seen as the stabiliser, in the affine
Weyl group $\Wh$, of the fundamental alcove (see \cite[Page
16]{Lam-Shi} or \cite[Page 5]{Magyar}). By
composition with the natural map $\Wh\to W$
to the finite Weyl group of $G$, one realises $Z$ as a subgroup of
$W$. This subgroup is given by the elements $(v_i)_{i\in I_m}$,  where
$v_i$ is the smallest element of $W$ such that
$v_i\varpi_i^\vee=w_0\varpi_i^\vee$,  where $w_0$ denotes
the longest element in $W$ (see \cite[Page 16]{Lam-Shi} for
example).

In a different context, Seidel proved in \cite{seidel} that the
fundamental group of the group of Hamiltonian symplectomorphisms of
a symplectic variety $X$ can be mapped to the group
$QH^*(X)_{loc}^\times$ of invertible elements of the quantum
cohomology ring localised in the quantum parameters (Seidel's
construction has coefficients in $\Z/2\Z$; here we use the
definition in \cite{mcduff} which has integer coefficients).
When $X=G/P$ is a rational homogeneous space, the natural projective
structure of $X$ induces a symplectic structure preserved by $G^{ad}$
(recall that the center $Z$ of $G$ acts trivially on $X$), and we get a map
$$Z=\pi_1(G^{ad})\longrightarrow QH^*(G/P)_{loc}^\times,$$
which we call Seidel's representation.
According to \cite{mcduff}, it is hard to describe this
map explicitly. This has been done by A. Postnikov for Grassmannians and the
variety of complete flags in type $A$ (see \cite{Po1} and \cite{Po2}).
In \cite{cmp}, we described explicitly Seidel's representation
when $X=G/P$ is a minuscule or cominuscule homogeneous space; in particular
we proved it is faithful. In this note we extend this result to all
homogeneous spaces. Moreover we provide an explicit formula for the product
by a class in the image of Seidel's representation, in the more
general setting of the equivariant quantum cohomology ring
$QH_T^*(X)_{loc}$.

To the  parabolic subgroup $P$ of $G$, containing a fixed Borel subgroup,
we associate the set  $I_P$ of vertices of the Dynkin diagram defining $P$
(with the convention that if $P$ is the Borel subgroup itself, $I_P=I$
is the whole set of vertices of the Dynkin diagram). For $w\in W$, we denote by
$\sigma^P(w)$ the Schubert class induced by $w$ in
$H^{2l(w)}(G/P,\Z)$ and by $\eta_P$ the natural surjection
$Q^\vee\to Q^\vee/Q^\vee_P$ from the coroot lattice to its quotient by
the coroot lattice of $P$. We prove:

\begin{theo}\label{main}
For any $i\in I_m$ and any $w\in W$, we have in $QH_T^*(G/P)_{loc}$:
$$\sigma^P(v_{i})* \sigma^P(w)=
q_{\eta_P(\varpi_{i}^\vee-w^{-1}(\varpi_{i}^\vee))}
\sigma^P(v_{i}w).$$
\end{theo}

\begin{theo}\label{seidel}
Let $f$ denote the Weyl involution. Seidel's representation of the
group $\pi_1(G^{ad})$ in $QH^*(G/P)_{loc}^\times$ is given, for
$i\in I_m$, by
$$-\varpi^\vee_i\mapsto  \sigma^P(v_{f(i)}).$$
\end{theo}

For the complete flag variety $G/B$, Theorem \ref{main} follows rather
directly from a factorization theorem in the affine Grassmannian due to
Magyar, and Peterson's comparison theorem between the equivariant homology
of the affine Grassmannian, and the quantum cohomology of $G/B$
(see Remark 3.18). In order to prove Theorem \ref{main} for arbitrary
flag varieties $G/P$, we will need to extend Magyar's result in a suitable
way. This extension will rely on the introduction of certain variants of the
tools introduced in \cite{Peterson} and \cite{Lam-Shi}, and an important
part of our work will consist in checking that some of the key statements
in these papers still hold in our extended setting. Once our generalized
Magyar's formula, Proposition \ref{magyar-gen}, is established, the proof of
Theorem 1 in 3.5 readily follows.

\section{Peterson's map and Magyar's factorisation formula}

\subsection{Affine Grassmannian, affine Weyl group and extended
  affine Weyl group}

Let us denote by $\Omega K$ the affine Grassmannian associated to
the group $G$. This is an ind-variety which can be defined as
follows. Let $K$ be a maximal compact subgroup of $G$; as a set, the
affine Grassmannian is the set of functions $f:S^1\to K$, such that
$f(1)=1_K$ and $f$ can be extended to a meromorphic function $f:D\to
G$ on the closed unit disk with poles only at the origin (see for
example \cite{Magyar} for more details). The affine Grassmannian
$\Omega K$ may also be seen as the quotient $G(\C((t)))/G(\C[[t]])$.

The affine Weyl group $\Wh$ is defined as the semidirect product of
the Weyl group $W$ by the coroot lattice $Q^\vee$  (see
\cite{Magyar} of \cite{Lam-Shi} for more details). In the same way
one defines the extended affine Weyl group $\Wt$ as the semidirect
product of $W$ by $P^\vee$ the coweight lattice (see \cite{Magyar}
or \cite{Lam-Shi}). Elements of these groups will be denoted by
$wt_\lambda$ with $w\in W$ and $\lambda\in P^\vee$ or $\lambda\in
Q^\vee$. We will denote by $P^\vee_-$ the set of antidominant
coweights and by $\Qt$ its intersection with the coroot lattice.

Note that $\Wh$ is naturally a subgroup of $\Wt$, which is normal,
and that the quotient $\Wt/\Wh$ is isomorphic to $P^\vee/Q^\vee$ and
thus to $Z$. The stabiliser of the fundamental alcove defines a
natural section $Z\to\Wt$ of this quotient (see \cite[Page
5]{Magyar}). The opposite of the minuscule coweights
$(-\varpi_i^\vee)_{i\in I_m}$ are representatives of the quotient
$P^\vee/Q^\vee$ in $P^\vee$. The image of $-\varpi_i^\vee$ with
$i\in I_m$ by this section is the element (see \cite[Page
16]{Lam-Shi}):
$$\tau_i:=v_it_{-\varpi_i^\vee}.$$

Recall the definition of the affine root lattice $Q_{\rm
  aff}=\bigoplus_{i=0}^n\Z\alpha_i$ where $\alpha_i$ are the simple
roots associated to the extended (or affine) Dynkin diagram . Recall
also the definition of the imaginary root $\delta=\alpha_0+\theta$
where $\theta$ is the longest root of the finite root system. For
more details on the affine root lattice, see \cite{Lam-Shi}. Recall
also the action of the extended Weyl group $\Wt$ on the affine root
lattice given by $wt_\lt(\alpha+n\delta)
=w(\alpha)+(n-\sca{\lt}{\alpha})\delta$. It follows that any
$\tau\in Z$ induces an automorphism $i\mapsto\tau(i)$ of the Dynkin
diagram defined by the formula $\tau(\alpha_i)=\alpha_{\tau(i)}$ for
any simple root $\alpha_i$.

We already defined the element $v_i$; it is also given by
$v_i=w_0w_0^{P_i}$, where $w_0$ denotes the longest element of the
Weyl group $W$, and $w_0^{P_i}$ the longest element of the Weyl
group $W_{P_i}$ of the maximal parabolic subgroup $P_i$ of $G$
associated to the cominuscule simple root $\alpha_i$\footnote{Recall
that these elements
  $w_0^{P_i}$ define by
  $\sigma^{P_i}(u)\mapsto\sigma^{P_i}(w_0^{P_i}u)$ the strange duality
  in $QH^*_{loc}(G/P_i)$, cf. \cite{cmp} and \cite{cmp2}}.
The element $v_i$ is also the longest element
in $W^{P_i}$, the set of minimal length representatives of the quotient
$W/W_{P_i}$. This element is well understood, in particular, we shall
use the following fact. Recall that the Weyl involution $f$ is defined
on roots by $\alpha\mapsto-w_0(\alpha)$ and induces an involution on
the set of simple roots, that we also denote by $f$.

\begin{fact}
\label{fait-utile}
  (\i) We have $v_i^{-1}=v_{f(i)}$.

(\i\i) Let $\alpha$ be a positive root, then $v_i(\alpha)$ is positive
if and only if $\sca{\varpi^\vee_i}{\alpha}=0$.
\end{fact}

%\noindent
%Recall the definition of the affine root lattice $Q_{\rm
%  aff}=\bigoplus_{i=0}^n\Z\alpha_i$ where $\alpha_i$ are the simple
%roots associated to the extended (or affine) Dynkin diagram . Recall
%also the definition of the imaginary root $\delta=\alpha_0+\theta$
%where $\theta$ is the longest root of the finite root system. For
%more details on the affine root lattice, see \cite{Lam-Shi}. From
%the previous fact and the action of the extended Weyl group $\Wt$ on
%the affine root lattice given by $wt_\lt(\alpha+n\delta)
%=w(\alpha)+(n-\sca{\lt}{\alpha})\delta$, we deduce the following
%useful fact.
%
%\begin{fact}
%  We have the formula $\tau(\alpha_i)=\alpha_{\tau(i)}$ for any simple
%  root $\alpha_i$.
%\end{fact}

\subsection{Equivariant homology}

We will be interested in the equivariant homology ring $H_*^T(\Omega
K)$ of $\Omega K$ for a maximal torus $T$ in $G$. We refer to
\cite{Magyar} for details. Let $S=H^*_T({\rm pt})$, then
$H_*^T(\Omega K)$ is a free $S$-module with basis $(\xi_x)_{x\in
\Waffm}$ where $\xi_x$ is the class of a Schubert variety in $\Omega
K$ and $\Waffm$ is the set of minimal length representatives of
$\Wh/W$ in $\Wh$. Note that we will assume, following \cite[Page
4]{Lam}, that $T$ does not contain the rotation action of $S^1$ on
$\Omega K$ (considered as a set of loops). This implies that the
image  in $S$ of  the imaginary root $\delta=\theta+\alpha_0$ is zero.

Let us denote by $R$ the root system, by $R^+$ the set of positive roots and
by $R_P^+$ the set of positive roots in $P$. Define
the subset $\Wpaff$ of $\Waff$ by
$$\Wpaff=\left\{wt_\lt\in\Waff\ /\ \forall\alpha\in R_P^+,\
  \sca{\lt}{\alpha}=\left\{ \begin{array}{ll}
0&\textrm{if } w(\alpha)>0\\
-1&\textrm{if } w(\alpha)<0\\
  \end{array}
\right.\right\}.$$
Recall that we have
$$
\Waff^- = \left\{ wt_\lambda \in \Waff \ /\ \forall i ,\
\sca{\lambda}{\alpha_i} = 0\Longrightarrow
 w(\alpha_i) > 0\right\}.
$$
Lam and Shimozono proved that
$$J_P=\bigoplus_{x\in\Waff\setminus\Wpaff}S\xi_x$$
is an ideal in $H_*^T(\Omega K)$. Moreover, they defined a map
$\pi_P:\Waff\to\Wpaff$ by proving that for any $x\in\Waff$, there
exists a unique factorisation $x=x_1x_2$ with $x_1\in\Wpaff$ and
$x_2\in(W_P)_{\rm aff}$ (the affine group associated to $W_P$). The
map $\pi_P$ is defined by $x\mapsto x_1$. We have the following
proposition (see \cite[Proposition 10.11]{Lam-Shi}):

\begin{prop}
\label{local-p}
  Let $x\in\Waffm\cap  \Wpaff$ and $\lambda\in\Qt$. Then
  $x\pi_P(t_\lt)\in\Waffm\cap \Wpaff$ and we have:
$$\xi_x\xi_{\pi_P(t_\lambda)}=\xi_{x\pi_P(t_\lambda)}\ \ {\rm mod}\ J_P.$$
\end{prop}

In particular, we may quotient and localise the equivariant homology
ring as follows, and define
$$H_*^T(\Omega K)_P:=(H_*^T(\Omega
K)/J_P)[\xi_{\pi_P(t_\lambda)}^{-1},\lambda\in\Qt].$$

\subsection{Equivariant quantum cohomology}

On the other hand, we consider the equivariant quantum cohomology
ring $QH^*_T(G/P)$. For details we refer to \cite{mihal} and
\cite{Lam-Shi}. This ring is a free $S[q_i,i\in I_P]$-module with
basis given by the Schubert classes $(\sigma^P(w))_{w\in W^P}$
(where $W^P$ is the set of minimal coset representatives of $W/W_P$
and $W_P$ is the Weyl group of $P$). We localise this ring as
follows:
$$QH^*_T(G/P)_{loc}:=QH^*_T(G/P)[q_i^{-1},i\in I_P].$$
Lam and Shimozono \cite[Theorem 10.16]{Lam-Shi} proved the so called
Peterson theorem (see also \cite{Peterson}):

\begin{theo}\label{qcoh}
  The map $\psi_P:H_*^T(\Omega K)_P\to QH^*_T(G/P)_{loc}$ defined by
  $$\xi_{w\pi_P(t_\lambda)}\xi_{\pi_P(t_\mu)}^{-1}\mapsto
  q_{\eta_P(\lambda-\mu)}\sigma^P(w),$$
  for $w\in W^P$ and $\lambda,\mu\in\Qt$, is an isomorphism of
  $S$-algebras.
\end{theo}

\subsection{Magyar's factorisation formula}

Denote by $\Wtaffm$ the set of minimal length
representatives of $\Wt/W$ in $\Wt$. An element $x\in\Wtaffm$ can be
uniquely written as $x=\tau_i\xh$ where $\xh\in\Waffm$ and $i\in I_m$.
Magyar proved the following formula (see \cite[Theorem
A]{Magyar}, the length condition $l(x)+l(t_\lt)=l(xt_\lt)$ being always
satisfied for $x$ in $\Waffm$ and $\lt\in P^\vee_-$):

\begin{theo}
  \label{magyar}
Let
$\lambda$ be antidominant in $P^\vee$ and let $x$ and $y$ in
$\Wtaffm$.
Assume that $xt_\lambda=y$. Then
$$\xi_{\xh}\xi_{\widehat{t_\lt}}=\xi_{\yh}\qquad
in \;\; H_*^T(\Omega K).$$
\end{theo}

In order  to prove Theorem \ref{main} we will need a
more general formula that will be proven in the next section. However
Magyar's formula was a guide for our proof, and it is sufficient to establish
 Theorem \ref{main} for  $G/B$ (see Remark \ref{ppmagyar}).

\section{Symmetries in the quantum cohomology}

Instead of  Magyar's formula, which is true in the full
equivariant homology ring $H^T_*(\Omega K)$, we prove a more general
formula in the quotient ring $H^T_*(\Omega K)/J_P$.
Theorem \ref{main} will be a straightforward application.

We start with the definition of an extended nil Hecke ring $\At$,
acting on the equivariant homology ring  $\Ht$ of the extended affine
Grassmannian $\Omega K^{\ad}$ (see section \ref{omegaKad}).
This action extends the known action of $\A$ on $\H$.

\subsection{Extended nil Hecke ring}

Because of the natural inclusion of $Z$ in $\Wt$, we can let $Z$ act
on the weight lattice; the result of the action of $\tau \in Z$ on a
weight $\lambda$ will be simply denoted $\tau(\lambda)$. We have
already observed that
$\tau(\alpha_i) = \alpha_{\tau(i)}$. Recall
 from \cite[Section 6.1]{Lam-Shi} the definition of the affine nil Hecke
ring. It is a non-commutative ring with generators $A_i$ for
$i \in \Ih$ and $\lambda$ in the weight lattice, subject to the
following relations:
$$
\begin{array}{rcll}
\lambda \mu & = & \mu \lambda\\
A_i \lambda & = & (s_i . \lambda)A_i + \langle \lambda , \alpha_i^\vee
\rangle \cdot 1\\
A_i A_i & = & 0\\
A_iA_jA_i \cdots & = & A_jA_iA_j \cdots & \mbox{ if }
s_is_js_i \cdots  =  s_js_is_j\cdots
\end{array}
$$
(in the last line, there are the same numbers of factors in all the
products). Note that the above set of relations is invariant under
the action of $\tau\in Z$, so that $\tau$ yields an algebra
automorphism of $\A$.

\begin{defi}
  The extended nil Hecke ring $\At$ is the smashed product
  $Z\ltimes \A$ of $\A$ by $Z$. As a $\Z$-module,
  it is just the tensor product $Z\otimes\A$. The
  multiplication is defined by
  $$(\tau\otimes a)(\sigma\otimes b)=\tau\sigma\otimes\sigma^{-1}(a)b.$$
\end{defi}

In other words, $\At$ is generated by $Z$ and $\Ah$, with the
relations $\tau a = \tau(a) \tau$.

\begin{prop}
  Let $x\in\Wt$ and write $x=\tau\xh$ with $\tau\in Z$ and
  $\xh\in\Waff$. Then the map $x\mapsto\tau\otimes\xh$ extends the
  inclusion of $\Waff$ in $\A$ to a multiplicative inclusion $\Wt\to\At$. The
  elements $x\in\At$ for $x\in\Wt$ form a basis (over ${\rm Frac}(S)$)
  of $\At$.
\end{prop}

\begin{proo}
  We only need to prove that this map is multiplicative. For
  $x=\tau\xh$ and $y=\sigma\yh$, we have
  $xy=\tau\sigma\sigma^{-1}(\xh)\yh$ in $\Wt$ and the result follows.
Since the elements of $\A$ defined by $\xh\in\Waff$ form a basis
of $\A$ over ${\rm  Frac}(S)$, we can deduce the result for $\At$. \end{proo}

\subsection{Equivariant homology of $\Omega K^\ad$}
\label{omegaKad}

Recall from \cite[Page 17]{Magyar} the definition of the variety
$\Omega K^{\ad}$ of based loops with values in the adjoint group
$K^{\ad}=K/Z$. The variety $\Omega K^{\ad}$ has its connected
components indexed by $Z$. There is an action of $T$ on $\Omega
K^\ad$ and the fixed points for this action are indexed by $\Wt/W$.
Define the space $LK^\ad$ of all loops in $K^\ad=K/Z$. We may also
(as for $\Omega K$, see \cite{Magyar}) see $\Omega K^\ad$ as the
quotient $LK^\ad/K^\ad$ by writing any loop $f\in LK^\ad$ as
$f(t)=f(t)f(1)^{-1}f(1)\in\Omega K^\ad\cdot K^\ad$. In particular,
we have a left action of $LK^\ad$ on $\Omega K^\ad$. The group $Z$
can be realised as the subgroup of $LK^\ad$ consisting of the points
defined by the loops
$(t\mapsto v_i\exp(-2i\pi t\varpi_i^\vee))$ (see
\cite[Page 6]{Magyar}). This together with the classical realisation
of $\Waff$ in $LK^\ad$ gives a realisation of the extended Weyl group
$\Wt$ in $LK^\ad$. The corresponding elements are precisely
the fixed points of  the $T$-action.

The equivariant cohomology ring $H^T_*(\Omega K^{\ad})$ has a basis
indexed by $\Wtaffm$: if $\xi_{\xh}$, with $\xh\in\Waffm$, is a class in
$H^T_*(\Omega K)$ and $\tau$ is an element in $Z$, then the
translation by $\tau\in LK^\ad$ in $\Omega K^{\ad}$ defines a class
$\xit_{\tau \xh}=\tau\cdot\xi_\xh$ in $H^T_*(\Omega K^\ad)$. Note that
$\xi_\tau$ is simply the class $[\tau]$. Following
\cite[Paragraph 1.2]{Magyar}, we have for any $\tau\in Z$ and any
$\xh\in\Waffm$:
$$\xi_\tau\xi_\xh=\xi_\xh\xi_\tau.$$
The ring $\Ht$ is thus the tensor product ring $Z\otimes\H$: the
product is defined by $$(\tau\otimes\xi_\xh)(\sigma\otimes\xi_\yh)=
\tau\sigma\otimes\xi_\xh\xi_\yh.$$

We may define operators $A_i$ acting on the ring $\Ht$, in the same way
as they are usually defined as operators on  $\H$.
That is, consider $E\to B$ the universal principal $LK^\ad$-bundle. Then
$LK^\ad$ acts on the right on $E\times \Omega K^\ad$ by
$f\cdot(e,f')=(ef^{-1},ff')$. In particular, $T$ acts on $E\times
\Omega K^\ad$ and we may consider the quotient $E\times^T \Omega
K^\ad$ of $E\times \Omega K^\ad$ by $T$. The equivariant cohomology
is defined by $H^*_T(\Omega K^\ad)=H^*(E\times^T \Omega K^\ad)$ and
the equivariant homology is the dual of that space. Now for any
vertex $i\in I\cup\{0\}$ of the extended Dynkin diagram, there
exists an associated compact subgroup $K_i^\ad$ of $LK^\ad$ (see
for example \cite[Page 15]{Magyar} for the subgroup $K_i$ of $LK$
the group of loops in $K$; the group $K_i^\ad$ is the image of $K_i$
in $LK^\ad$). Consider the map $\pi_i:E\times^T \Omega K^\ad\to
E\times^{K^\ad_i} \Omega K^\ad$ and define the operator $A_i$ by integrating
on the fibers and pulling-back.  By
duality this defines an action of $A_i$ on $\Ht$.
We may also consider  any $s\in S$ as an operator on $\Ht$. It is enough
to do this for a weight $\lt$, for which the action is defined  by
intersection with $c_1({\cal L}_\lt)$,
where ${\cal L}_\lt$ denotes the line bundle defined by $\lt$. Finally we
define the operator $\tau$ for any $\tau\in Z$ thanks to the
translation by $\tau$ in $\Omega K^\ad$.

\begin{prop}
These operators define an action of $\At$ on $\Ht$ extending the
action of $\A$ on $\H$. This action can be written as follows:
$$(\tau\otimes a)\cdot(\sigma\otimes\xi)=
\tau\sigma\otimes(\sigma^{-1}(a)\cdot\xi).$$
\end{prop}

\begin{proo}
If the action is well defined, it certainly extends the action of $\A$ on
$\H$. We only need to verify that the commutation relations $\tau
A_i\tau^{-1}=A_{\tau(i)}$ and $\tau s\tau^{-1}=\tau(s)$ are
satisfied.

For the first one, we only need to remark that the conjugate group
$\tau(K_i)\tau^{-1}$ is $K_{\tau(i)}$. For the second one, we may suppose that
$s$ is a  weight $\lt$. Let ${\cal L}'_\lt$ be the line bundle defined
by the twisted action of $T$ on $\C$ given for $t\in T$ by
multiplication with $\tau(\lt)(t)$. The operator $\tau s\tau^{-1}$ is
given by the intersection with $c_1({\cal L}'_\lt)$. This
coincides with the intersection with $c_1({\cal L}_{\tau(\lt)})$.
\end{proo}

For $x\in\Wt$, write $x=\tau\xh$ with $\tau\in Z$ and $\xh\in\Waff$.
We define the element $\at_x\in\At$ by $\at_x=\tau\otimes A_\xh$.
We also define the length of an element $x = \tau \xh$ as the number
of inverted positive real roots, as in \cite[Section 10.1]{Lam-Shi},
so that $l(x) = l(\xh)$.

\begin{prop}
\label{action}
  The action of $\At$ on $\Ht$ is given by the following formula
$$\at_x\cdot\xit_y=\left\{
\begin{array}{cl}
  \xit_{xy}&\textrm{if }l(xy)=l(x)+l(y)\ {\rm and}\ xy\in\Wtaffm\\
0&\textrm{otherwise.}
\end{array}\right.$$
\end{prop}

\begin{proo}
  Write $x=\tau\xh$ and $y=\sigma\yh$ with $\tau,\sigma\in Z$ and
  $\xh,\yh\in\Waff$. We have $xy=\tau\sigma\sigma^{-1}(\xh)\yh$. We
  compute
$$
\begin{array}{cl}
  \at_x\cdot\xit_y&=(\tau\otimes A_\xh)\cdot(\sigma\otimes\xi_y)\\
&=\tau\sigma\otimes \sigma^{-1}(A_\xh)\cdot\xi_y\\
&=\tau\sigma\otimes A_{\sigma^{-1}(\xh)}\cdot\xi_y.\\
\end{array}$$
But we know by a result of Kostant and Kumar \cite{KK} (see also
\cite[Section 6.2]{Lam-Shi}) that
$$A_{\sigma^{-1}(\xh)}\cdot\xi_\yh=\left\{
\begin{array}{cl}
  \xi_{\sigma^{-1}(\xh)\yh}&\textrm{if
  }l(\sigma^{-1}(\xh)\yh)=l(\sigma^{-1}(\xh))+l(\yh)\ {\rm and}\
  \sigma^{-1}(\xh)\yh\in\Waffm\\
0&\textrm{otherwise,}
\end{array}\right.$$
thus we have the equality
$$\at_x\cdot\xit_y=\left\{
\begin{array}{cl}
  \xit_{xy}&\textrm{if }l(\sigma^{-1}(\xh)\yh)=l(\xh)+l(\yh)\ {\rm
    and}\ \sigma^{-1}(\xh)\yh\in\Waffm\\
0&\textrm{otherwise.}
\end{array}\right.$$
The equality $l(\sigma^{-1}(\xh))=l(\xh)$ comes from the fact that
$\sigma$ permutes the simple roots and the formula
$l(\xh)=\vert\{\alpha\in R_{\rm aff}^+\ /\ \xh(\alpha)<0\}\vert$.
Moreover, since $xy=\tau\sigma\sigma^{-1}(\xh)\yh$, we know that $xy$
is in $\Wtaffm$ if and only if $\sigma^{-1}(\xh)\yh$ is in $\Waffm$.
Furthermore, $l(\sigma^{-1}(\xh)\yh)=l(xy)$ and
$l(x)+l(y)=l(\xh)+l(\yh)$ so that the condition
$l(\sigma^{-1}(\xh)\yh)=l(\xh)+l(\yh)$ is equivalent to the
condition $l(xy)=l(x)+l(y)$.
\end{proo}

\begin{rema}
There exist, a priori, two actions of the extended affine Weyl
group $\Wt$ on $H_*^T(\Omega K^\ad)$: one is given by the embedding
of $\Wt$ in $\At$. The other one is simply defined
on $H_*^T(\Omega K^\ad)$ by the operators $R_x^*$ induced, for $x\in\Wt$,
by the right multiplication $R_x$ in $E\times^T \Omega K^\ad$, by $x$
considered as an element of  $\subset LK^\ad$. This right multiplication
operator is well defined
because $\Wt$ normalises $T$. It is a classical (but non trivial)
result that these two actions coincide for the classical affine Weyl
group $\Waff$ (see \cite[Theorem 11.3.9]{kumar}). Since the action of
$\tau$ on a class $\xi_x$ is given by $\xi_{\tau x}$ (see
Proposition \ref{action}), which is also the class obtained by
right translation on $E^\ad\times^T \Omega K^\ad$, these two actions
of $\Wt$ coincide as well.
\end{rema}

\subsection{The map $\jad$}

We want to generalise Theorem 6.2 in \cite{Lam-Shi}. For this we
need to consider another basis of the ring $\Ht$ given by $T$-fixed
points in $\Omega K^\ad$. Recall that these fixed points are indexed
by $\Wt/W$ and that good representatives are given by $P^\vee$ by
$$p_{t_\lt}=(t\mapsto\exp(2i\pi t\lt))\qquad \mathrm{for}\;\;
\lt\in P^\vee .$$
Furthermore, for $\lt$ and $\mu$ in $P^\vee$, these fixed points
satisfy $p_{t_\lt} p_{t_\mu}=p_{t_{\lt}t_{\mu}}$ and the affine Weyl
group action gives $x\cdot p_{t_\lt}=p_{xt_\lt}$.

\begin{defi}
Let $\psit_{t_\lambda}$ for $\lt\in P^\vee$ be the element in $\Ht$
defined by
$$\psit_{t_\lambda}= i_\lt^* : H^*_T(\Omega K^{ad}) \rightarrow S,$$
where $i_\lt$ denotes the inclusion of the $T$-fixed
point $p_\lt$ in $\Omega K^\ad$.
\end{defi}

The fact that the elements $\psit_{t_\lt}$, for $\lt\in P^\vee$, form a
basis of $\Ht$ over ${\rm Frac}(S)$, comes from the same
statement for compact spaces (see \cite[Theorem C8, page 537]{kumar}),
and the fact that the extended affine Grassmannian, as well as the
affine Grassmannian, is an increasing union of compact $T$-stable
finite dimensional Hausdorff subspaces.

\begin{prop}
\label{mult-psit} We have the following formulas:
$$\psit_{t_\lt}\psit_{t_\mu}=\psit_{t_{\lt+\mu}}\ \ \ \textrm{and}
\ \ \ x\cdot\psit_{t_{\lt}}=\psit_{xt_{\lt}}$$ for any $\lt$ and
$\mu$ in $P^\vee$ and any $x$ in $\Wt$.
\end{prop}

\begin{proo}
These formulas are immediate consequences of the identities  $p_{t_\lt}
p_{t_\mu}=p_{t_{\lt}t_{\mu}}$ and $x\cdot p_{t_\lt}=p_{xt_\lt}$.
\end{proo}

We are now in position to generalise \cite[Theorem 6.2]{Lam-Shi}. We
will denote by $Z_\At(S)$ the centraliser of $S$ in $\At$. We also
define the ideal $\Jt$ in $\At$ by
$$\Jt=\sum_{w\in W\setminus\{{\rm id}\}}\At(1\otimes A_w).$$

\begin{prop}
\label{6.2-mod}
  There is an $S$-algebra isomorphism $\jad:\Ht\to Z_\At(S)$ such that
  for any $x$ and $y$ in $\Wtaffm$ we have $\jad(\xit_x)=\at_x\ {\rm
    mod}\ \Jt$ and $\jad(\xit_x)\cdot\xit_y=\xit_x\xit_y$.
\end{prop}

\begin{proo}
This extends the similar result
in \cite{Lam}, whose proof we follow.
We define $\jad$ by letting $\jad(\psit_{t_\lt})=t_\lt$ for $\lt\in
P^\vee$, and extending by ${\rm Frac}(S)$-linearity. Proposition
\ref{mult-psit} shows
that it is an $S$-algebra morphism. The formula
$x\cdot\psit_t=\psit_{xt}$ gives
$$\jad(\psit_{t_\lt})\cdot\psit_{t_\mu}=\psit_{t_\lt t_\mu}
=\psit_{t_\lt}\psit_{t_\mu}$$
which, since the elements $\psit_{t_\lt}$
for $\lt\in P^\vee$ form a basis of $\Ht$, implies
the identity $\jad(\xit_x)\cdot\xit_y=\xit_x\xit_y$.

To prove that the image is contained in $Z_\At(S)$, consider the
action
$t_\lt\cdot(\alpha+n\delta)=(\alpha+n\delta-\sca{\alpha}{\lt}\delta)t_\lt$
where $\alpha$ is a root and
$\delta=\alpha_0+\theta$ is the imaginary root. As  the image of
$\delta$ in $S$ is zero, the commutation relation follows.

We have $\jad(\xit_x)\cdot\xit_\id=\xit_x\xit_\id=\xit_x$ for
$x\in\Wtaffm$. This implies that $\jad(\xit_x)=\at_x+a$ where $a$ lies
in the annihilator of $\xit_\id\in\Ht$ in $\At$.
This annihilator is clearly the ideal $\Jt$.

Finally, let $a=\sum_{x\in\Wt}a_xx$ be an element in $Z_\At(S)$,
then we have for $v\in S$
$$a v=\sum_{x\in\Wt}a_xx(v)x=va=\sum_{x\in\Wt}a_xvx.$$
In particular $a_x=0$ for all $x$ such that $x(v)\neq v$ for some
$v\in S$. But recall that the action of $w\in W$ is given by
$w(\alpha+n\delta)=w(\alpha)+n\delta$ and because $\delta$ is send
to zero, this equals $\alpha$ if and only if $w=\id$. In particular
$x$ has to be a translation, which exactly means it belongs to
the image of $\jad$.
\end{proo}

\subsection{Generalised Magyar's formula}

In this subsection, we generalise Magyar's formula to the quotient
$H^T_*(\Omega K)/J_P$. This formula, as for Magyar's formula, is a
generalisation to coweights of the localisation formula of
Proposition \ref{local-p}.
Let us first introduce some definitions.

\begin{defi}
    We define the subset $\Wtpaff$ of $\Wt$ as follows~: for
    $x=\tau\xh$ with $\tau\in Z$ and $\xh\in\Waff$, $x\in\Wtpaff$
    if $\xh\in\Wpaff$.
\end{defi}

\begin{lemm}
Define $(R_P)_{\rm
  aff}^+=\{\alpha+n\delta\in R^+_{\rm aff}\ /\ \alpha\in R_P\}$, then
we have
$$\Wtpaff=\{x\in\Wt\ /\ x(\beta)>0\ \textrm{for all}\
\beta\in(R_P)_{\rm aff}^+\}$$
\end{lemm}

\begin{proo}
Recall from \cite[Section 10.3]{Lam-Shi}, that the condition in the
lemma defines $\Wpaff$ in $\Waff$. The extension to $\Wtpaff$
 follows from the fact
that $Z$ stabilises the cone of positive roots.
\end{proo}

For an element $x\in\Wt$, we write $x=\tau\xh$ with $\tau\in Z$ and
$\xh\in\Waff$. We need to compare the properties of $x$ and $\xh$. We
first compare their expressions as products $x=wt_\lt$ with
$\lt\in P^\vee$, and
$\xh=vt_\mu$ with $\mu\in Q^\vee$.

\begin{lemm}
\label{p-r}
  Let $w\in W$ and $\lt\in P^\vee$. Write
  $\lt=-\varpi_{i(\lt)}^\vee+\lh$, where $\lh\in Q^\vee$ and $i(\lt)\in
  I_m$. Then
$$wt_\lt=\tau_{i(\lt)}v_{f(i(\lt))}wt_{\lt-w^{-1}(\varpi_{f(i(\lt))}^\vee)},$$
where $\lt-w^{-1}(\varpi_{f(i(\lt))}^\vee)$ is in $Q^\vee$.
\end{lemm}

\begin{proo}
  We know that there exist an index $j$, an element $u\in W$  and an
  element $x\in Q^\vee$ such that $wt_\lt=\tau_jut_x$. We then get
  the equalities $u=v_j^{-1}w$ and $\mu=\lt-w^{-1}(\varpi_{f(j)})$. But
  $\mu$ has to be in $Q^\vee$, hence $f(j)=f(i(\lt))$ and
  thus $j=i(\lt)$.
\end{proo}

Let $\lambda\in P^\vee$ and $w\in W$. Write
$wt_\lambda=\tau_{i(\lambda)}v_{f(i(\lambda))}w
t_{\lambda-w^{-1}(\varpi_{f(i(\lambda))}^\vee)}$ with $\lh\in
Q^\vee$.

\begin{lemm}
For any positive root $\alpha\in R^+$ of the finite root system, we
have the equality
$$\chi(w(\alpha)<0)+\sca{\lt}{\alpha} =
\chi(v_{f(i(\lt))}w(\alpha)<0)+\sca{\lt-w^{-1}
  (\varpi^\vee_{f(i(\lt))})}{\alpha},$$
where $\chi(A)=1$ if $A$ is true and $\chi(A)=0$ otherwise.
\end{lemm}

\begin{proo}
Let us consider the following four cases:
\begin{itemize}
\item[1.] $w(\alpha)>0$ and $\sca{\varpi^\vee_{f(i(\lt))}}{w(\alpha)}=1$;
\item[2.] $w(\alpha)>0$ and $\sca{\varpi^\vee_{f(i(\lt))}}{w(\alpha)}=0$;
\item[3.] $w(\alpha)<0$ and $\sca{\varpi^\vee_{f(i(\lt))}}{w(\alpha)}=-1$;
\item[4.] $w(\alpha)<0$ and $\sca{\varpi^\vee_{f(i(\lt))}}{w(\alpha)}=0$.
\end{itemize}
According to these cases and setting $u=v_{f(i(\lambda))}w$ and
$\mu=\lambda-w^{-1}(\varpi_{f(i(\lambda))}^\vee)$, we have:
\begin{itemize}
\item[1.] $u(\alpha)<0$ and $\sca{\mu}{\alpha}=\sca{\lt}{\alpha}-1$;
\item[2.] $u(\alpha)>0$ and $\sca{\mu}{\alpha}=\sca{\lt}{\alpha}$;
\item[3.] $u(\alpha)>0$ and $\sca{\mu}{\alpha}=\sca{\lt}{\alpha}+1$;
\item[4.] $u(\alpha)<0$ and $\sca{\mu}{\alpha}=\sca{\lt}{\alpha}$.
\end{itemize}
The statement can then be checked case by case.
\end{proo}

This lemma gives a generalisation to the extended affine Weyl group
$\Wt$ of the length formula given in \cite[Lemma 3.1]{Lam-Shi}:

\begin{coro}
\label{gen-long}
  Let $wt_\lt\in\Wt$ with $w\in W$ and $\lt\in P^\vee$. Then
$$l(wt_\lt)=\sum_{\alpha\in R^+}\vert\chi(w(\alpha)<0)+\sca{\lt}{\alpha}\vert.$$
\end{coro}

\begin{proo}
  The formula $wt_\lambda=\tau_{i(\lambda)}v_{f(i(\lambda))}w
t_{\lambda-w^{-1}(\varpi_{f(i(\lambda))}^\vee)}$ implies that
$$l(wt_\lt)=l(v_{f(i(\lambda))}w
t_{\lambda-w^{-1}(\varpi_{f(i(\lambda))}^\vee)}),$$ and \cite[Lemma
3.1]{Lam-Shi} gives the formula
$$l(v_{f(i(\lambda))}wt_{\lambda-w^{-1}(\varpi_{f(i(\lambda))}^\vee)})=\sum_{\alpha\in
  R^+}\vert\chi(v_{f(i(\lambda))}w(\alpha)<0)+
\sca{{\lambda-w^{-1}(\varpi_{f(i(\lambda))}^\vee)}}{\alpha}\vert.$$
We conclude thanks to the previous lemma.
\end{proo}

The previous lemma also provides a criterion, for an element
$x\in\Wt$, to be such that $\xh\in\Wpaff$. More precisely let
$wt_\lt\in\Wt$. Write once more
$wt_\lambda=\tau_{i(\lambda)}v_{f(i(\lambda))}w
t_{\lambda-w^{-1}(\varpi_{f(i(\lambda))}^\vee)}$.

\begin{coro}
\label{wpaff-p}
  The element $v_{f(i(\lambda))}w
  t_{\lambda-w^{-1}(\varpi_{f(i(\lambda))}^\vee)}$ is in $\Wpaff$ if
  and only if, for any root $\alpha\in R^+_P$,
$$\sca{\lt}{\alpha}=\left\{
  \begin{array}{ll}
0&\textrm{if } w(\alpha)>0\\
-1&\textrm{if } w(\alpha)<0.\\
  \end{array}
\right.$$
\end{coro}

\begin{proo}
This condition is exactly the condition for an element
$wt_\lt\in\Waff$ to be in $\Wpaff$, see Lemma 10.1 in
\cite{Lam-Shi}. But remark that this condition is equivalent to
$$\chi(w(\alpha)<0)+\sca{\lt}{\alpha}=0$$
for all $\alpha\in R^+_P$. We conclude once more by the
previous lemma.
\end{proo}

\begin{coro}
An element $wt_\lt$ with $w\in W$ and $\lt\in P^\vee$ is in $\Wtpaff$
if and only if,  for any root $\alpha\in R^+_P$,
$$\
  \sca{\lt}{\alpha}=\left\{ \begin{array}{ll}
0&\textrm{if } w(\alpha)>0\\
-1&\textrm{if } w(\alpha)<0.\\
  \end{array}
\right.$$
\end{coro}

\begin{proo}
This is an immediate application of Corollary \ref{wpaff-p}.
\end{proo}

The map $\pi_P:\Waff\to\Wpaff$ extends to $\pit_P:\Wt\to\Wtpaff$ by
setting $\pi_P(x)=\tau\pi_P(\xh)$ for $x=\tau\xh$.
Many results of \cite[Paragraph 10]{Lam-Shi} on $\Wpaff$ extend readily to
$\Wtpaff$. This is the case of Lemma 3.1 (see Corollary
\ref{gen-long}), Lemmas 3.2 and 3.3 and Lemmas 10.2 to 10.6.
The definition of the extension $\pit_P$ of $\pi_P$ ensures that
the results of Proposition 10.8 in \cite{Lam-Shi} remain true in our setting.

Furthermore, we extend the definition of the ideal $J_P$ by defining
the following ideal $\Jt_P$ of $\At$:
$$\Jt_P=\sum_{x\in\Wtaffm\setminus\Wtpaff}S\xit_x.$$
This is indeed an ideal because $\Wtaffm$ and $\Wtpaff$ are stable under left
multiplication by $Z$.

\begin{prop}
\label{magyar-gen}
  For all $x\in\Wtaffm\cap  \Wtpaff$ and $\lambda\in P^\vee_-$, we have
  $x\pi_P(t_\lt)\in\Wtaffm\cap \Wtpaff$ and
$$\xi_\xh\xi_{\widehat{\pi_P(t_\lambda)}}=
\xi_{\widehat{x\pi_P(t_\lambda)}}\ \ {\rm mod}\ J_P.$$
\end{prop}

\begin{proo}
  We follow the proof of Proposition \ref{local-p} by Lam and
  Shimozono (see \cite[Proposition 10.11]{Lam-Shi}).
Lemma 10.10 in \cite{Lam-Shi} extends readily to

  \begin{lemm}
    Let $\lt\in P^\vee$ be antidominant. Then
    $(1\otimes A_i)\cdot\xit_{\pi_P(t_\lt)}=0\ {\rm mod}\ \Jt_P$ for all
    $i\in I$.
  \end{lemm}

\noindent In particular $\Jt\cdot \xi_{{\pi_P(t_\lambda)}}=0\ {\rm
mod}\ \Jt_P$. Now we can apply Proposition \ref{6.2-mod} to get the identity:
$$\xit_x\xit_{{\pi_P(t_\lambda)}}=\at_{x}\cdot \xit_{{\pi_P(t_\lambda)}}
\ \;{\rm mod}\ \Jt_P.$$ Let us prove that the product $x\pi_P(t_\lt)$ is
length additive. We have seen that
Proposition 10.8 of \cite{Lam-Shi} extends to $\Wt$
and $\Wtpaff$. In particular $x=w\pi_P(t_\nu)$ for $w\in W^P$ and
$\nu\in P^\vee$ antidominant and we have
$x\pi_P(t_\lt)=w\pi_P(t_\nu)\pi_P(t_\lt)=w\pi_P(t_{\nu+\lt})$. By
\cite[Lemma 3.3]{Lam-Shi}, because $x$ and $x\pi_P(t_\lt)$ are in
$\Wtaffm$, we have $l(x)=l(\pi_P(t_\nu))-l(w)$ and
$l(x\pi_P(t_\lambda))=l(\pi_P(t_{\nu+\lt}))-l(w)$. We only need to
prove that $l(\pi_P(t_{\nu+\lt}))=l(\pi_P(t_\nu))+l(\pi_P(_\lt))$.
Because Lemmas 3.3, 10.3 and 10.6 of \cite{Lam-Shi} extend to
coweights, the same proof as for coroots in \cite{Lam-Shi} gives the
additivity. By proposition \ref{action}, we thus have the formula
$\xit_x\xit_{\pi_P(t_\lt)}=\xit_{x\pi_P(t_\lt)}\ {\rm mod}\ \Jt_P$. But writing
$x=\tau\xh$, $\pi_P(t_\lt)=\sigma\widehat{\pi_P(t_\lt)}$ and
$x\pi_P(t_\lt)=\tau\sigma\widehat{x\pi_P(t_\lt)}$, this implies the
equality
$(\tau\otimes\xi_\xh)(\sigma\otimes\xi_{\widehat{\pi_P(t_\lt)}})=
(\tau\sigma\otimes\xi_{\widehat{x\pi_P(t_\lt)}})\ {\rm mod}\ \Jt_P$
and thus
$\xi_\xh\xi_{\widehat{\pi_P(t_\lt)}}=
\xi_{\widehat{x\pi_P(t_\lt)}}\ {\rm mod}\ J_P$.

Finally let us prove that $x\pi_P(t_\lt)$ belongs to
$\Wtaffm\cap \Wtpaff$. Since
$x\in\Wtpaff$, we have $x\pi_P(t_\lt)=\pi_P(x)\pi_P(t_\lt)=\pi_P(xt_\lt)$
(from the generalisation of \cite[Proposition 10.8]{Lam-Shi} to $\Wt$). But
$x\in\Wtaffm$ and $\lt\in P^\vee_-$, thus, by the generalisation of
\cite[Lemma 3.3]{Lam-Shi} to $\Wt$, we get that
$xt_\lt\in\Wtaffm$. Thanks, once more,  to the generalisation of
\cite[Proposition 10.8]{Lam-Shi} to $\Wt$, we can then deduce the result.
\end{proo}

\subsection{Application to symmetries}

We are now in position to prove
Theorem \ref{main} for any $G/P$. Let $w$
be an element of the Weyl group $W$, let $\varpi_i^\vee$ be a
minuscule coweight and let $\mu$ and $\nu$ be in $Q^\vee$. We begin with
the formulas:
$$t_{-\varpi_i^\vee-\mu}=\tau_iv_{f(i)}
t_{-(\varpi_i^\vee+\varpi_{f(i)}^\vee+\mu)}\ {\rm and}\
wt_{-\nu}t_{-\varpi_i^\vee-\mu}=\tau_iv_{f(i)}w
t_{-(\varpi_i^\vee+w^{-1}(\varpi_{f(i)}^\vee)+\mu+\nu)}.$$ Applying
the map $\pi_P$, we get the equalities
$\pi_P(t_{-\varpi_i^\vee-\mu})=\tau_i\pi_P(v_{f(i)})
\pi_P(t_{-(\varpi_i^\vee+\varpi_{f(i)}^\vee+\mu)})$ and
$\pi_P(wt_{-\nu}t_{-\varpi_i^\vee-\mu})=\tau_i\pi_P(v_{f(i)}w)
\pi_P(t_{-(\varpi_i^\vee+w^{-1}(\varpi_{f(i)}^\vee)+\mu+\nu)})$. For
$\mu$ and $\nu$ dominant enough, the elements $wt_{-\nu}$, $v_{f(i)}
t_{-(\varpi_i^\vee+\varpi_{f(i)}^\vee+\mu)}$ and $v_{f(i)}w
t_{-(\varpi_i^\vee+w^{-1}(\varpi_{f(i)}^\vee)+\mu+\nu)}$ all belong to
$\Waffm$, and their images by $\pi_P$ are in $\Wpaff\cap \Waffm$. We
may thus apply Proposition \ref{magyar-gen}, and we obtain the following
formula:
$$\xi_{\pi_P(w)\pi_P(t_{-\nu})}
\xi_{\pi_P(v_{f(i)})\pi_P(t_{-(\varpi_i^\vee+\varpi_{f(i)}^\vee+\mu)})}=
\xi_{\pi_P(v_{f(i)}w)\pi_P(t_{-(\varpi_i^\vee+w^{-1}(\varpi_{f(i)}^\vee+\mu+\nu))})}\ \;
{\rm mod}\ J_P.$$ Applying Peterson's map of Theorem
\ref{qcoh}, we get the corresponding formula in the quantum cohomology ring:
$$\sigma^P(w)q_{-\eta_P(\nu)}*\sigma^P(v_{f(i)})
q_{-\eta_P((\varpi_i^\vee+\varpi_{f(i)}^\vee+\mu))} =
\sigma^P(v_{f(i)}w)q_{-\eta_P((\varpi_i^\vee+w^{-1}(\varpi_{f(i)}^\vee)
  +\mu+\nu))},$$
hence finally:
$$\sigma^P(w)*\sigma^P(v_{f(i)}) =
q_{\eta_P(\varpi_{f(i)}^\vee-w^{-1}(\varpi_{f(i)}^\vee))}
\sigma^P(v_{f(i)}w).$$ This concludes the proof of Theorem \ref{main}.

\begin{rema}
\label{ppmagyar}
 As a final remark in this section, let us explain briefly
why Magyar's formula is not sufficient to prove Theorem \ref{main}
in full generality. The point is that, in order to
apply this formula, we need an element $t_\lt\in\Wtaffm\cap\Wtpaff$.
But Corollary \ref{wpaff-p} implies the following fact:

\begin{fact}
A translation $t_\lt$, with $\lt\in P^\vee$, is in $\Wtaffm\cap\Wtpaff$
if and only if $\lt\in P^\vee_-$ and $\sca{\lt}{\alpha}=0$ for all
$\alpha\in R_P^+$.
\end{fact}

In particular, if $P$ is a maximal parabolic subgroup associated to
a cominuscule coweight $\varpi_i^\vee$ (equivalently, if $I_P=\{i\}$
for some $i$ in $I_m$), this implies that $\lt$ needs to be a multiple
of $\varpi_i^\vee$.  But then Magyar's formula gives the Theorem \ref{main}
only for multiples of $v_{f(i)}$. The corresponding symmetries of the
quantum cohomology ring are those generated by the quantum product
with the punctual  class $[{\rm pt}]$, which were described in
\cite{cmp}. These symmetries  do not generated the full group $Z$ in
general.

However, as we already mentionned, if the parabolic subgroup $P$ is
a Borel subgroup, there is no restriction on $\lt\in P^\vee_-$. In this case
our Theorem \ref{main} follows from Magyar's formula.
\end{rema}

\section{Seidel's representation}

In this section, we prove Theorem \ref{seidel}. That is, we prove
that $\sigma^P(v_{f(i)})$ is the invertible element
$S(\omega_i^\vee)$ of the (localized) quantum cohomology ring,
corresponding to $-\omega_i^\vee \in Z=\pi_1(G^\ad)$ through
Seidel's representation. As we recalled in the introduction, we use
the integer-valued definition given in \cite[Example 8.6.8]{mcduff}.
Let us  recall Seidel's construction.

Since $G$ is simply-connected, the coweight lattice parametrise
1-parameter subgroups of $T$, our prefered maximal torus of $G$
(which is contained in $P$). Let $\C^* \subset T$ correspond to
$\omega_i^\vee$. The image of $S^1 \subset \C^* \subset T \subset G
\rightarrow G^{ad}$ defines an element of $\pi_1(G^{ad})$ which will be
denoted $\pi_1(\omega_i^\vee)$.

We now set $M_i = (\C^2 - \{0\}) \times^{\C^*} G/P$. There is a natural map
$\pi : M_i \rightarrow (\C^2 - \{0\})/\C^* \simeq \p^1$.
The fibers of $\pi$ are isomorphic with $G/P$. As in \cite{seidel}, we denote
$TM_i^v$ the vertical tangent space, that is, the kernel of the differential
$d\pi$ in $TM_i$. Choose a point $z \in \p^1$, and if
$\cal S$ is a space of sections of $\pi$, denote by
$ev:{\cal S} \rightarrow \pi^{-1}(z) \simeq G/P$ the
evaluation map at $z$. In the following proposition, $k$ denotes the number
of roots $\alpha$ in the unipotent radical of $P$ such that
$\scal{\omega_i^\vee,\alpha} > 0$.

\begin{prop}    \label{prop_sections_degre_minimal}
For any section  $s$ of $\pi$, one has $\deg (TM_i^v)(s) \geq -k$.
Moreover, let $\cal S$ denote the space of sections $s$ of $\pi$ such that
$\deg(TM_i^v)(s)=-k$. Then
$ev_* [{\cal S}] \in H_*(G/P)$ is the class $\sigma^P(v_{f(i)})$ of a
homogeneous Schubert variety.
\end{prop}

We postpone the proof  to the end of the
section. Let us first investigate the geometry of $M_i$ in more details. For
$\lambda$ a $W_P$-invariant element of the root lattice, we denote
by ${\cal L}_\lambda$ the associated $G$-linearised bundle on $G/P$.
Moreover, if ${\cal L}$ is any $G$-linearised bundle, we denote
${\cal L}^i$ the bundle ${\cal L} \times^{\C^*} (\C^2-0)$ over
$M_i$. Let us finally denote $\lambda_P = \sum_{\alpha \in \Phi^+ -
\Phi_P} \alpha$.

\begin{prop}
We have $\det (TM_i^v) = {\cal L}_{\lambda_P}^i$.
\end{prop}
\begin{proo}
Let $\eta : G/P \times (\C^2-0) \rightarrow M_i$ denote the quotient map,
and let $p_1 : G/P \times (\C^2-0) \rightarrow G/P$ denote the first
projection. Since $\eta$ is a submersion, it induces an
isomorphism of vector bundles $p_1^* TG/P \simeq \eta^* TM_i^v$. Moreover, this
isomorphism commutes with the $\C^*$-action, so that it induces an
isomorphism $TG/P \times^{\C^*} (\C^2-0) \simeq TM_i^v$.
We deduce  that
$\det(TM_i^v) \simeq \det(TG/P) \times^{\C^*} (\C^2-0) \simeq
{\cal L}_{\lambda_P}^i$.
\end{proo}

\begin{coro}
For $s:\p^1 \rightarrow M_i$ any section of $\pi$, one has
$\deg s^* \det(TM_i^v) \geq -k$.
\end{coro}

\begin{proo}
Let ${\cal N} = \det(TM_i^v) \otimes \pi^* \cO(k)$: we will show that
$\cN$ is nef. This implies that for any map $f:\p^1 \rightarrow M_i$,
$\deg\ f^* \cN \geq 0$. Hence for any section $s$ of $\pi$ we deduce
that $\deg\ s^* \det(TM_i^v) \geq -k$.

Let $\lambda_1,\ldots,\lambda_u$ be the
weights of the $G$-module
$\Gamma ( G/P , \det(TG/P) )$, counted with multiplicities
(so that $u = \dim \Gamma ( G/P , \det(TG/P) )$. Since
$\pi_* \det(TM_i^v) = \Gamma ( G/P , \det(TG/P) ) \times^{\C^*} (\C^2-0)$,
we get
\begin{equation}
\label{pi*N}
\pi_* \cN = \bigoplus_{j=1}^r
\cO(\scal{\lambda_j,\omega_i^\vee} + k).
\end{equation}

For the lowest weight $-\lambda_P$ of $\Gamma ( G/P , \det(TG/P) )$,
we have $\scal{-\lambda_P , \omega_i^\vee} = - k$. Therefore all the
line bundles in (\ref{pi*N}) have non-negative degree. But then
$\cN$ is base-point free, hence nef. This concludes the proof.
\end{proo}

Since $\det(TM_i^v)$ is relatively ample, we get an embedding $M_i
\rightarrow \p (\pi_* \cN)^\vee$. Let us decompose $\Gamma(G/P , \det(TG/P) )$
as $\bigoplus_{j\geq 0} E_j$,
where our $\C^*$ has weight $j-k$ on $E_j$. We have an inclusion
$$(G/P
\cap \p E_0^\vee) \times \p^1 \subset M_i \subset \p (\pi_*
\cN)^\vee.$$

\begin{prop}
If $s$ is a section of $\pi$
such that $\deg\ s^*\ \det(TM_i^v) = -k$, there exists
$x_0 \in G/P \cap \p E_0^\vee$
such that $s(u) = (u,x_0) \in \p^1 \times (G/P \cap \p E_0^\vee)$ for all
$u \in \p^1$.
\end{prop}

\begin{proo}
The space of sections of $\cN$ identifies with $\bigoplus E_j
\otimes S^j \C^2$. Let $f : M_i \rightarrow \p (\bigoplus E_j^\vee
\otimes S^j \C^2)$ be the morphism defined by $\cN$. Since $f^*
\cO(1) = \cN$, $\deg\ s^* \det(TM_i^v) = -k$ implies that $f \circ s$
is constant.

Note that for $u \not = v \in \p^1$, we have
$f(\pi^{-1}(u)) \cap f(\pi^{-1}(v)) \subset \p E_0^\vee$. Since $f \circ s$
is constant, this implies that for all $u \in \p^1$, $s(u)$ belongs to
$\p^1 \times \p E_0^\vee$ and that the induced morphism
$\p^1 \rightarrow \p E_0^\vee$ is constant. The proposition is
proved.
\end{proo}

Proposition \ref{prop_sections_degre_minimal} is a consequence of
the last two results and the following lemma \ref{E0} (in fact, we
still have to prove that $\cal S$ is reduced; this will be a
consequence of the following argument).

Let us compute the constant term of $S(\omega_i^\vee)$, the Seidel
element corresponding to the loop $\omega_i^\vee$. With the
notations of \cite{mcduff}, this is $S_0(\pi)$ where $\pi : M_i
\rightarrow \p^1$ is our fibration. According to the remark before
\cite[Proposition 7.11]{seidel}, we can compute $S_0(\pi)$ as the
push-forward $(ev_{z_0})_* {\cal S}$ of the space ${\cal S}$ of
sections of $\pi$ of degree $-k$; by proposition
\ref{prop_sections_degre_minimal} we deduce $S_0(\pi) =
d.\sigma^P(v_{f(i)})$, where $d \geq 1$, with $d>1$ in case $\cal S$
is not reduced.

We can now easily complete the proof of Theorem \ref{seidel}.
We can write any class $x \in QH^*(G/P)$ as  $x = \sum x_{w,I}
\sigma^P(w)q^I$ for some integers $x_{w,I}$. We  let $|x| := \sum
x_{w,I}$. By formula \cite[8.6.4]{mcduff} and since $M_i$ is
an algebraic variety, all the coefficients ${S(\omega_i^\vee)}_{w,I}$ are
non-negative. Moreover we know that ${S(\omega_i^\vee)}_{v_{f(i)},0} =
d$. From Theorem \ref{main}, and since the product of any two
effective classes is again effective, it follows that for any
effective class $x$ in $QH^*(G/P)$, we have $|x*S(\omega_i^\vee)| \geq
|x|$. Since $S(\omega_i^\vee)$ is unipotent, this forces
$|S(\omega_i^\vee)| = 1$, hence  $S(\omega_i^\vee) =
\sigma^P(v_{f(i)})$. In particular $d=1$, which also completes the proof of
proposition \ref{prop_sections_degre_minimal}.

\smallskip All that remains to prove is the following claim:

\begin{lemm}   \label{E0}
Let ${\cal L}$ be a very ample homogeneous line bundle on $G/P$. Let
$G/P \subset \p \Gamma(G/P,{\cal L})^\vee$ be the corresponding
embedding. Let $E_0 \subset \Gamma(G/P,{\cal L})$ be the lowest
weight space for a one parameter subgroup of $G$ corresponding to a
simple root $\alpha_i$, and $P_i \subset G$ the corresponding
parabolic subgroup. Then $(G/P \cap \p E_0^\vee)_{red}$ is the
closed $P_i$-orbit in $G/P$; its homology class is
$\sigma^P(v_{f(i)})$.
\end{lemm}

\begin{proo}
Let $x \subset \Gamma(G/P,{\cal L})^\vee$ be the highest weight
line. Then $x \in G/P$. Moreover, $\Gamma(G/P,{\cal L})^\vee$ is
generated as a vector space by $U(\g).x$, where $U(\g)$ denotes the
enveloping algebra of $\g$. Consider $\displaystyle \llie = \tlie \oplus
\bigoplus_{\alpha:\langle \alpha,\omega_i^\vee \rangle = 0 }
\g_\alpha$, a Levi subalgebra of $\plie_i=Lie(P_i)$; let $L \subset
G$ be the corresponding subgroup. By definition of $E_0$, its dual $E_0^\vee$
is  generated by the lines $U(\llie).x$. This implies
that $E_0^\vee$ is an irreducible $L$-module.

The closed $L$-orbit $L.x$ is contained in $G/P \cap \p E_0^\vee$.
We claim that they are equal. In fact, let ${\cal O}$ be any
$L$-orbit in $G/P \cap \p E_0^\vee$. We have $L.x \subset
\overline{\cal O}$. Moreover,
$$T_x L.x  = \bigoplus_{\langle \alpha,
\omega_i^\vee \rangle = 0} \g_\alpha . x= T_x \overline{\cal O}.$$
Thus $L.x= \overline{\cal O}$.

Finally, in order to identify the Schubert class $[L.x]$, we consider the
incidence diagram
$$
\begin{array}{ccc}
G/(P\cap P_i) & \rightarrow & G/P\\
\downarrow p_i\\
G/P_i
\end{array}
$$
Let $w \in W_{P_i}$ be the element corresponding to a fiber of $p_i$
(that is, the minimum length representative of the class modulo
$W_P$ of the longest element in $W_{P_i}$), and let $v^{P \cap P_i}
\in W$ represent the open orbit in $G/(P \cap P_i)$. The element in
$W^{P_i}$ corresponding to the open orbit in $G/P_i$ is $v_i$, and
therefore we have the relation $v^{P \cap P_i} = v_i w$, which
implies, since $v_i^{-1} = v_{f(i)}$, the relation $w = v_{f(i)}
v^{P \cap P_i}$. Let $F_i$ denote a fiber of $p_i$ in $G/(P \cap
P_i)$. We therefore have, by definition, $[F] = \sigma^{P \cap
P_i}(v_{f(i)})$. Since the projection $G/(P \cap P_i) \rightarrow
G/P$ restricts to an isomorphism on $F_i$, we deduce that $[L.x] =
\sigma^P(v_{f(i)})$.
\end{proo}

\bigskip\noindent
Pierre-Emmanuel {\sc Chaput}, Laboratoire de Math{\'e}matiques Jean Leray,
UMR 6629 du CNRS, UFR Sciences et Techniques,  2 rue de la Houssini{\`e}re, BP
92208, 44322 Nantes cedex 03, France.

\noindent {\it email}: pierre-emmanuel.chaput@math.univ-nantes.fr

\medskip\noindent
Laurent {\sc Manivel},
Institut Fourier, UMR 5582 du CNRS,  Universit{\'e} de Grenoble I,
BP 74, 38402 Saint-Martin d'H{\`e}res, France.

\noindent {\it email}: Laurent.Manivel@ujf-grenoble.fr

\medskip\noindent
Nicolas {\sc Perrin},  Institut de Math{\'e}matiques, Universit{\'e}
Pierre et Marie Curie, Case 247, 4 place Jussieu, 75252 PARIS Cedex
05, France and Hausdorff Center for Mathematics, Universit{\"a}t Bonn,
Landwirtschaftskammer (Neubau) Endenicher Allee 60, 53115 Bonn,
Germany.

\noindent {\it email}: nperrin@math.jussieu.fr

\end{document}